\documentclass{amsart}
\usepackage{amssymb, amsmath}
\usepackage[dvips]{graphicx}
\usepackage{amsfonts}
\usepackage{latexsym}
\usepackage{color}
\usepackage{pifont}
\newtheorem{theorem}{Theorem}	
\newtheorem{lemma}{Lemma}[section]		
\newtheorem{corollary}{Corollary}[section]		
\newtheorem{proposition}{Proposition}[section]

\title[]
{
Constant diameter and constant width of spherical convex bodies
}
\author{
Huhe Han}
\address{College of Science, Northwest Agriculture and Forestry University, China}
\email{han-huhe@nwafu.edu.cn}
\author{Denghui Wu}
\address{College of Science, Northwest Agriculture and Forestry University, China}
\email{wudenghui66@163.com}
\begin{document}
\begin{abstract}
In this paper we show that a spherical convex body $C$ is of constant diameter $\tau$ if and only if $C$ is of constant width $\tau$, for $0<\tau<\pi$. Moreover, some applications to Wulff shapes are given.
\end{abstract}
\maketitle
\section{Introduction}
Let $S^{n}$ denote the unit sphere of the $(n+1)$-dimensional Eucliean space $\mathbb{R}^{n+1}$.
For a convex body $K$ in $\mathbb{R}^{n+1}$,
it is well-known that the property that $K$ is of constant diameter is equivalent to
that it is of constant width. It is natural to ask if this fact
holds in the other spaces. This note considers the question in the unit sphere $S^{n}$.

A set $K\subseteq S^n$ is said to be {\it spherically convex}, if
$\{\lambda u:\lambda>0, u\in K\}$ is convex in $\mathbb{R}^{n+1}$.
If the body $K$ is closed, spherically convex and has nonempty interior
(with respect to $S^n$), then the body $K$ is spherical convex body.
The boundary of $K$ is denoted by $\partial K$.

For an integer $1\leq k\leq n$, $H$ denotes a $(k+1)$-dimensional subspace of
$\mathbb{R}^{n+1}$.
We call the intersection of the unit sphere $S^{n}$ with $H$ to be
a $k$-{\it dimensional subsphere} of $S^{n}$.

For two points $P,Q\in S^{n}$, the natural spherical distance of $P$ and $Q$
can be given by
\[
|PQ|=\arccos(\overrightarrow{OP}\cdot\overrightarrow{OQ}),
\]
where $O$ is the origin, and $\cdot$ denotes
the standard Euclidean scalar product.
Let $K$ be a spherical convex body in $S^{n}$. The \emph{diameter} of convex
body $K$ is defined by $\max\{|PQ|: P,Q\in K\}$.
Following \cite[Part 4]{LM18}, we say a convex body $K$ in the sphere
is of \emph{constant diameter} $\tau$,
if the diameter of $K$ is $\tau$, and
for every point $P\in\partial K$ there exists a point $Q$ of $K$ such that
$|PQ|=\tau$.

Let  $P$ be a point of $S^{n}$. For the $P$, we denote by $S_P$ the hypersphere with centered at $P$,
which is given by
\[
S_P=\{Q\in S^{n}:\overrightarrow{OP}\cdot\overrightarrow{OQ}=0\}.
\]
The closed hemisphere with centered at $P$ is written by $S^+_P$,
which is
\[
S^+_P=\{Q\in S^{n}:\overrightarrow{OP}\cdot\overrightarrow{OQ}\geq0\}.
\]
For any subset $W$ of $S^n$, the set
\[
\bigcap_{P\in W}S^+_P
\]
is called the spherical polar set of $W$ and is denoted by $W^{\circ}$.
For any non-empty closed hemispherical subset $W\subset S^n$, the
equality $\mbox{s-conv}(W)=(\mbox{s-conv}(W))^{\circ\circ}$ holds,
see, e.g. \cite{nishimurasakemi2}.
Here  $\mbox{s-conv}(W)$ is the spherical convex hull of $W$:
\[
\mbox{s-conv}(W)=\left\{\frac{\sum_{i=1}^{k}t_iP_i}{||\sum_{i=1}^{k}t_iP_i||_2}:
\sum_{i=1}^{k}t_i=1,t_i\geq 0, P_i\in W, k\in \mathbb{N}\right\},
\]
where $||\cdot||_2$ denotes the standard $n$-dimensional Euclidean norm.
For more details about spherical polar sets, see e.g. \cite{hnams,hnjmsj,hnsm,nishimurasakemi2}.
Let $K$ be a spherical convex body, and let $P$ be a point on the boundary of $K$.
If the hemisphere $S^+_{Q}$ contains $K$ and $P\in\partial K\cap S^+_{Q}$,
we say $S^+_{Q}$ {\it supports} $K$ {\it at} $P$, or that $S^+_{Q}$ is a {\it supporting
hemisphere} of $K$ at $P$. If for every point $P\in\partial K$ there exists a
unique hemisphere supporting $K$, then $K$ is called smooth.

If hemispheres $S^+_{P}$ and $S^+_{Q}$ of $S^n$ are different and
not opposite, then $S^+_{P}\cap S^+_{Q}$ is called a {\it lune} of $S^n$.
The {\it thickness} of lune $S^+_{P}\cap S^+_{Q}$ is given by
$\Delta(S^+_{P}\cap S^+_{Q})=\pi-|PQ|$.
If $S^+_{P}$ is a supporting hemisphere of a spherical convex body $K$,
the \emph{width} of $K$ with respect to $S^+_{P}$ is defined in \cite{Lassak15} by
\[
\text{width}_{S^+_{P}}(K)=\min\{\Delta(S^+_{P}\cap S^+_{Q}): K\subset S^+_{Q}\}.
\]
The thickness of $K$ is the minimum of $\text{width}_{S^+_{P}}(K)$
over all supporting hemispheres $S^+_{P}$.
We say that the spherical convex body $K$ is of \emph{constant width},
if all widths of $K$ with respect to any supporting hemispheres
are equal. The width of $W$ is denoted by $\Delta(W)$.

In this note we will show that

\begin{theorem}\label{maintheorem}
Let $C$ be a spherical convex body in $S^n$, and $0<\tau<\pi$.
The following two are equivalent:
\begin{enumerate}
	\item $C$ is of constant diameter $\tau$.
	\item $C$ is of constant width $\tau$.
\end{enumerate}
\end{theorem}
The Theorem \ref{maintheorem} was motivated by \cite{Lassak19}. We refer to the reader \cite{Lassak19} for the cases of smoothness boundary and $S^2$.
The statement from the proposition (2) to (1) has been given in \cite{LM18}.
We consider in this note for cases without smoothness of boundary and for $S^n$.

\par
\bigskip
This paper is organized as follows. In Section 2, the proof of Theorem \ref{maintheorem} is given, and some applications to Wulff shapes are given in Section 3.
\section{Proof of Theorem \ref{maintheorem}}
The following statement is one side part of Theorem \ref{maintheorem}
which was given in \cite{LM18}.
\begin{lemma}\cite[Theorem 5]{LM18}\label{diameter=a}
	If $C\subset S^n$ is a spherical convex body of constant width $\tau$, then
$C$ is of constant diameter $\tau$.
\end{lemma}

The next fact is useful in the coming proofs.
\begin{lemma}\label{lemma2.2}
Let $C$ be a spherical convex body in $S^n$.
For any point $\widetilde{P}$ in $\partial C^{\circ}$,
the center of the supporting hemisphere of $C^{\circ}$ at $\widetilde{P}$
is in $\partial C$.
\end{lemma}
\par
\begin{proof}
Let $\widetilde{P}$ be a point in $ \partial C^{\circ}$ and let
$S_P^+$ be the supporting hemisphere of $C$ at $\widetilde{P}$.
Since $C^\circ$ is a subset of $S_P^+$, it follows that
$\mid P\widetilde{R}\mid\leq\pi /2$, for all $\widetilde{R}\in C^\circ$.
By the fact $C^{\circ\circ}=C$, it follows that $P$ is a point in $C$ .
\par
Suppose that $P$ is an interior point of $C$.
Then there exists a sufficiently small real number $\varepsilon$
such that $\big(B(P,\varepsilon)\cap S^n\big)\subset C.$
Since $S_P^+$ is the supporting hemisphere of $C$ at $\widetilde{P}$, we have
$\mid P\widetilde{P}\mid=\pi /2$.
Without loss of generality, set $\widetilde{P}=(0,\dots,0,1)$ and $P=(1,0,\dots 0)$.
Since the point
\[Q=\left(\cos \frac{\varepsilon}{2}, 0,\dots, 0,-\sin \frac{\varepsilon}{2}\right)\]
satisfies $\mid PQ\mid=\arccos(\overrightarrow{OP}\cdot\overrightarrow{OQ})=\frac{\varepsilon}{2} <\varepsilon$, we have $Q\in\big(B(P,\varepsilon)\cap S^n\big)\subset C$.
From $\widetilde{P}\in\partial C^\circ$ and  $C=\bigcap_{\widetilde{R}\in C^\circ}S^+_{\widetilde{R}}$,
it follows $Q\in\bigcap_{\widetilde{R}\in C^\circ}S^+_{\widetilde{R}}\subset S^+_{\widetilde{P}}$.
However, from $\overrightarrow{O\widetilde{P}}\cdot\overrightarrow{OQ}<0$,
it follows that $Q$ is not contained in $S^+_{\widetilde{P}}$, which
is a contradiction.
Therefore, $P$ is a boundary point of $C$.

\end{proof}

Using Lemmas \ref{diameter=a} and \ref{lemma2.2}, we have
the next equivalent assertions for spherical convex bodies in $S^n$.
\begin{lemma}\label{widthofduals}
Let $C$ be a spherical convex body in $S^n$, and $0<\tau<\pi$.
The following two assertions are equivalent:
\begin{enumerate}
	\item $C$ is of constant width $\tau$.
	\item $C^{\circ}$ is of constant width $\pi-\tau$.
\end{enumerate}
\end{lemma}
\par
We refer to the reader \cite{michal} for the case of $S^2$.
\begin{proof}
Since $C=C^{\circ\circ}$ for any spherical convex body in $S^n$,
it is sufficient to prove that the assertion (1) implies the assertion (2).
We assume that $C$ is of constant width $\tau$.
For any point $\widetilde{P}$ in the boundary $\partial C^{\circ}$,
we let $S_{P}^+$
be a supporting hemisphere of $C^{\circ}$ at $\widetilde{P}$,
where $P$ is a point in the boundary of $C$ by Lemma \ref{lemma2.2}.
There is a supporting hemisphere $S_{Q}^+$ of $C^{\circ}$ with $Q\in \partial C$,
such that
\[
\mbox{width}_{S^+_{P}}(C^{\circ})=\Delta(S^+_{P}\cap S^+_Q)=\pi-|PQ|.
\]
Since $C$ is of constant width $\tau$, and $P,Q$ are two points in $\partial C$,
Theorem \ref{diameter=a} implies
$|PQ|\leq \tau.$
Then, we have
\[
\mbox{width}_{S^+_{P}}(C^{\circ})=\Delta(S^+_{P}\cap S^+_Q)=\pi-|PQ|\geq \pi-\tau.
\]
Suppose that $\mbox{width}_{S^+_{P}}(C^{\circ})>\pi-\tau$.
For any point $R$ in $\partial C$, $C^{\circ}$ is contained in
$S^+_R$. It follows that
\[
\pi-|PR|=\Delta(S^+_{P}\cap S^+_R)\geq\mbox{width}_{S^+_{P}}(C^{\circ})>\pi-\tau.
\]
This means the sharp inequalities
$|PR|<\tau$ hold for the fixed $P$ and any points $R$ in $\partial C$.
This contradicts to the fact that $C$ is of constant diameter $\tau$, which
follows from Lemma \ref{diameter=a} and our assumption that $C$ is of
constant width $\tau$.
Therefore, $\mbox{width}_{S^+_{P}}(C^{\circ})=\pi-a$.
\end{proof}

We now in the position to show our main theorem.

\begin{proof}
We firstly prove that the assertion (1) implies the assertion (2).
We assume that $C$ is of constant diamater $\tau$.
By Lemma \ref{widthofduals}, it is sufficient to show that
$C^{\circ}$ is a spherical convex body of constant width $\pi-a$.
Let $\widetilde{P}$ be a point in $\partial C^{\circ}$, here $\partial C^\circ$
stands for the boundary of the spherical polar body $C^\circ$.
Let $S^+_{P}$ be the supporting hemisphere of $C^{\circ}$ at $\widetilde{P}$.
Lemma \ref{lemma2.2} implies that $P$ is a point in $\partial C$.
By assumption, there exists a point $Q$ in $\partial C$ such that $|PQ|=\tau$.
It is clearly that $C^{\circ}\subset (S^+_{P}\cap S^+_{Q})$ and
\begin{equation*}\label{eq1}
\tag{$\star$} \mbox{width}_{S^+_{P}}(C^{\circ})\leq \Delta(S^+_{P}\cap S^+_Q)=\pi-|PQ|=\pi-a.
\end{equation*}
We claim that the equality in ($\star$) hold. In fact, there will be
a contradiction if $\mbox{width}_{S^+_{P}}(C^{\circ})<\pi-\tau$.
Suppose that $\mbox{width}_{S^+_{P}}(C^{\circ})=b<\pi-\tau$.
Then there is a supporting hemisphere $S^+_R$ of $C^{\circ}$ such that
\[
\mbox{width}_{S^+_{P}}(C^{\circ})=\Delta(S^+_{P}\cap S^+_R)=\pi-|PR|=b.
\]
This implies
\begin{equation*}\label{eq2}
\pi-|PR|=b=\mbox{width}_{S^+_{P}}(C^{\circ})<\pi-\tau,
\end{equation*}
or $|PR|>\tau$.
This contracts with that $\tau$ is the diameter of $C$,
since $P$, $R$ are points in $C$. Then we have $\mbox{width}_{S^+_{P}}(C^{\circ})=\pi-\tau$.
By the arbitrary selection of $P$, we get $C^\circ$ is
of constant width $\pi-\tau$.
By Lemma \ref{widthofduals}, we get that $C$ is of constant
width $\tau$.\\
\indent
Nextly, we use Lemma \ref{lemma2.2} to give a simple proof
of Lemma \ref{diameter=a}. We assume that $C$ is of constant width $\tau$.
Lemma \ref{widthofduals} shows that
$C^{\circ}$ is of constant width $\pi-\tau$.
Then, for every point $\widetilde{P}$ in $\partial C^{\circ}$,
there exists $P$ in $\partial C$ such that
$\mbox{width}_{S^+_{P}}(C^{\circ})=\pi-\tau$.
For the fixed $P$,
there exists a point $Q\in\partial C$ so that
$\mbox{width}_{S^+_{P}}(C^{\circ})=\Delta(S^+_{P}\cap S^+_Q)=\pi-\tau$,
which shows $|PQ|=\tau$.
Moreover,
$\pi-\tau=\mbox{width}_{S^+_{P}}(C^{\circ})\leq\Delta(S^+_{P}\cap S^+_R)=\pi-|PR|$,
for any point $R\in\partial C$.
Thus $|PR|\leq \tau$ for any point $R\in\partial C$.
Then the diameter of $C$ is $|PQ|=\tau$.
The result follows by the arbitrary selection of $P$.
\end{proof}

As a consequence of Theorem \ref{maintheorem} and Lemma \ref{widthofduals},
we have
\begin{corollary}\label{cor2.4}
Let $C$ be a spherical convex body in $S^n$, and $0<\tau<\pi$.
The following two propositions are equivalent:
\begin{enumerate}
	\item $C$ is of constant diameter $\tau$.
	\item $C^\circ$ is of constant diameter $\pi-\tau$.
\end{enumerate}
\end{corollary}

\section{Applications to Wulff shapes}
We first give some basic definitions related to Wulff shapes.

\subsection{Wulff shapes}
\indent
Let $\gamma:S^n \to \mathbb{R}_+$ be a continuous function, where $\mathbb{R}_+$ is the set of positive real numbers.
We denote by $H_{\theta,\gamma}^-$ the half space determined by $\gamma$ and $\theta\in S^n$,
\[
H_{\theta,\gamma}^-=\{x\in \mathbb{R}^{n+1}: x\cdot \theta\leq \gamma(\theta)\}.
\]
Then the {\it Wulff shape} associated with the function $\gamma$, denoted by $\mathcal{W}_\gamma$, is
defined by
\[
\bigcap_{\theta\in S^n} H_{\theta,\gamma}^-.
\]
Since the Wulff shape $\mathcal{W}_\gamma$ is a convex body in $\mathbb{R}^{n+1}$,
the radial function $\rho_{{}_{\mathcal{W}_\gamma}}$ of $\mathcal{W}_\gamma$ is given by
$\rho_{{}_{\mathcal{W}_\gamma}}(\theta)=\max\{\lambda>0:\lambda\theta\in\mathcal{W}_\gamma\},$
for $\theta\in S^n$.
Let $\overline{\gamma}:S^{n}\to \mathbb{R}_{+}$
: continuous function,
defined by
$\overline{\gamma}(\theta)=1/\rho_{{}_{\mathcal{W}_\gamma}}(-\theta)$.
The Wulff shape $\mathcal{W}_{\bar{\gamma}}$ associated with the support function
$\bar{\gamma}$ is called {\it dual Wulff shape of
$\mathcal{W}_\gamma$}, denoted by $\mathcal{DW}_\gamma$, namely, $\mathcal{W}_{\bar{\gamma}}=\mathcal{DW}_{\gamma}$.
We call $\mathcal{W}$ is a {\it self-dual Wulff shape} if $\mathcal{W}$ and its dual Wulff shape $\mathcal{DW}$ are exactly the same body, namely, $\mathcal{W}=\mathcal{DW}$.

Note that the dual Wulff shape $\mathcal{DW}_{\gamma}$ is related to
the Euclidean polar body of Wulff shape $\mathcal{W}_{\gamma}$ with
\begin{equation*}\label{eq4.1}
\tag{$\star\star$}\mathcal{DW}_{\gamma}=-(\mathcal{W}_{\gamma})^\circ,
\end{equation*}
where the Euclidean polar body $K^\circ$ of convex body $K$ in $\mathbb{R}^{n+1}$ is defined by
$$K^\circ=\{x\in\mathbb{R}^{n+1}:x\cdot y\leq1\}.$$
For more details on Wulff shapes, see for instance \cite{han1, hnams, Schneider14, Wu17}.

\subsection{Spherical Wulff shapes}

Let $Id: \mathbb{R}^{n+1}\to \mathbb{R}^{n+1}\times \{1\}\subset \mathbb{R}^{n+2}$
be the mapping defined by
$$Id(x)=(x,1).$$
Let $N$ denote the north pole of $S^{n+1}$, i.e., $N=(0,\ldots,0,1)\in \mathbb{R}^{n+2}$,
and let $S_{N,+}^{n+1}$ denote the north open hemisphere of $S^{n+1}$, i.e.,
$S_{N,+}^{n+1}=S_{N}^{+}\backslash S_{N}=\{Q\in S^{n+1}:\overrightarrow{ON}\cdot\overrightarrow{OQ}>0\}.$
The central projection relative to $N$, denoted by
$\alpha_N: S_{N,+}^{n+1}\to \mathbb{R}^{n+1}\times \{1\}$, is defined by
\[
\alpha_N\left(P_1, \ldots, P_{n+1}, P_{n+2}\right)
=
\left(\frac{P_1}{P_{n+2}}, \ldots, \frac{P_{n+1}}{P_{n+2}}, 1\right).
\]
We call the spherical convex body $\widetilde{W}_\gamma=\alpha^{-1}(Id(\mathcal{W_\gamma}))$ is {\it the spherical Wulff shape of $\mathcal{W_\gamma}$}.

The diameter and width of spherical Wulff shape and their polar bodies have
the following conclusion from Theorem \ref{maintheorem}, Lemma \ref{widthofduals} and
Corollary \ref{cor2.4}.
\begin{corollary}
Let $\gamma:S^n\to \mathbb{R}_+$ be a continuous function. If the spherical Wulff shape $\widetilde{W}_\gamma=\alpha_N^{-1} \circ Id(\mathcal{W_\gamma})$ of $\mathcal{W}_\gamma$ is of constant width. Then
\begin{enumerate}
\item $\Delta(\widetilde{W}_\gamma)\ +\mbox{diam}\ (\widetilde{W}_\gamma^\circ)=\pi$,
\item $\Delta(\widetilde{W}_\gamma)\ + \Delta\ (\widetilde{W}_\gamma^\circ)=\pi$,
\item $\mbox{diam}(\widetilde{W}_\gamma)\ +\Delta\ (\widetilde{W}_\gamma^\circ)=\pi$,
\item $\mbox{diam}(\widetilde{W}_\gamma)\ +\mbox{diam}\ (\widetilde{W}_\gamma^\circ)=\pi$,
\end{enumerate}
where $\Delta(A)$ and $\mbox{diam}(A)$ are the width and the diameter of spherical convex body $A$
in $S^n$, respectively.
\end{corollary}

The following proposition gives a relationship between Euclidean dual Wulff shapes
and spherical polar bodies of Wulff shapes.
\begin{proposition}\cite{nishimurasakemi2}\label{dualbysphere}
Let $\gamma:S^n\to \mathbb{R}_+$ be a continuous function. Then
\[
\mathcal{DW}_\gamma=Id^{-1}\circ \alpha_N \big((\alpha_N^{-1} \circ Id(\mathcal{W}_\gamma))^\circ \big).
\]
\end{proposition}

This together with \eqref{eq4.1} shows that
$$\alpha_N^{-1} \circ Id(-\mathcal{W}^\circ)=(\alpha_N^{-1} \circ Id(\mathcal{W}))^\circ ,$$
namely, the order of the map $\alpha_N^{-1} \circ Id$ and polar operation can be exchanged up to a sign.

By Proposition \ref{dualbysphere}, it follows that a Wulff shape $\mathcal{W}$ is self-dual if and only if $\alpha_N^{-1} \circ Id(\mathcal{W}) =(\alpha_N^{-1} \circ Id(\mathcal{W}))^\circ $. Moreover,
\begin{proposition}\cite{hnjmsj}\label{selfdual}
Let $\gamma:S^n\to \mathbb{R}_+$ be a continuous function.
Then $\mathcal{W}_\gamma$ is a self-dual Wulff shape if and only
if its spherical Wulff shape is of constant width $\pi/2$, namely,
the spherical convex body $\alpha_N^{-1} \circ Id(\mathcal{W_\gamma})$
is of constant width $\pi/2$.
\end{proposition}

This combining with Theorem \ref{maintheorem} implies:

\begin{corollary}
Let $\gamma:S^n\to \mathbb{R}_+$ be a continuous function.
Then $\mathcal{W}_\gamma$ is a self-dual Wulff shape if and only if
its spherical Wulff shape is of constant diameter $\pi/2$, namely,
the spherical convex body $\alpha_N^{-1} \circ Id(\mathcal{W_\gamma})$ is of constant diameter $\pi/2$.
\end{corollary}

{\bf Acknowledgements.}
This work was supported, in partial, by
Natural Science Basic Research Plan in Shaanxi Province of China
(2019JQ-246),
and the Initial Foundation for Scientific Research
of Northwest A\&F University (2452018016, 2452018018).


\end{document}